\documentclass{article}

\usepackage{graphicx}
\usepackage{amscd}
\usepackage{amsmath}
\usepackage{amsthm}
\usepackage{amsfonts}
\usepackage{amssymb}

\usepackage{color}
\usepackage{ulem}

\newtheorem{theorem}{Theorem}

\newtheorem{corollary}{Corollary}

\newtheorem{example}{Example}

\newtheorem{lemma}{Lemma}

\newtheorem{proposition}{Proposition}

\newcommand{\A}{\mathcal{A}}
\newcommand{\B}{\mathcal{B}}
\newcommand{\C}{\mathcal{C}}

\newcommand{\F}{\mathbb{F}}
\newcommand{\M}{\mathbb{M}}
\newcommand{\N}{\mathbb{N}}
\newcommand{\T}{\mathbb{T}}

\DeclareMathOperator{\Gal}{Gal}

\DeclareMathOperator{\Aut}{Aut}
\DeclareMathOperator{\Out}{Out}
\DeclareMathOperator{\Inn}{Inn}

\title{Arbitrarily large Galois orbits of non-homeomorphic surfaces}

\author{Gabino Gonz\'{a}lez-Diez \and Gareth A. Jones \and David Torres-Teigell}

\begin{document}

\maketitle

\begin{abstract}
\noindent We construct orbits of the absolute Galois group, of explicit unbounded size, consisting of surfaces with mutually non-isomorphic fundamental groups. These are Beauville surfaces with Beauville group $PGL_2(p)$.

\medskip

(secondary).

\end{abstract}

\section{Introduction}

If $X$ is a projective variety defined over the field $\overline{\mathbb{Q}}$ of algebraic numbers, then it is natural to ask which topological properties it shares with its Galois conjugates $X^{\sigma}$, obtained by applying elements $\sigma$ of the absolute Galois group $\Gal\overline{\mathbb{Q}}/\mathbb{Q}$ to the coefficients of the polynomials defining $X$. In dimension~$1$ the answer is straightforward: two Galois conjugate curves have the same genus, and are therefore homeomorphic, so they share all their topological properties; for this and other Galois invariants, in the context of Grothendieck's theory of {\it dessins d'enfants}, see~\cite{JSt}.

More generally, by Hodge theory (see~\cite{GH,Voi} for instance) the dimensions of the cohomology groups $H^i(X,\mathbb{C})$ of a complex projective variety $X$ can be expressed in terms of the Hodge numbers $h^{p,q}(X) = \dim H^q(X,\Omega^p)$, where $\Omega^p$ is the sheaf of holomorphic $p$-differential forms on $X$. By Serre's GAGA principle~\cite{Ser56} these numbers $h^{p,q}$ are invariant under Galois conjugation. It follows that in dimension $2$, for instance, many of the standard topological invariants of a complex projective surface $X$ are also Galois invariants. These include
\begin{itemize}
\item the Betti numbers $b_i=\dim H_i(X)$,
\item the Euler characteristic or Euler number $e=\sum_{i=0}^4(-1)^ib_i$,
\item the irregularity $q=h^{0,1}=h^{1,0}$,
\item the geometric genus $p_g=h^{0,2}=h^{2,0}$,
\item the arithmetic genus $p_a=p_g-q$,
\item the holomorphic Euler characteristic $\chi=p_g-q+1$,
\item the signature (of the second cohomology group) $\tau=4\chi-e$, and
\item the Chern numbers $c_2=e$ and $c_1^2=K^2=12\chi-e$.
\end{itemize}
(see e.g.~\cite[Th.~6.33]{Voi}). Nevertheless, in 1964 Serre~\cite{Ser64} constructed examples of Galois conjugate pairs of complex projective varieties, including surfaces, which have non-isomorphic fundamental groups, and are therefore not homeomorphic to each other. Since then, further examples of conjugate but non-homeomorphic varieties have been found: see~\cite{Abe, ACC, EV, Cha, MS, BCGarX, Shi, Fue} for instance.

A Beauville surface is an example of a complex surface which is isogenous to a higher product, that is, it has the form $S = (C_1 \times C_2)/G$ where each $C_i$ is a curve of genus $g_i>1$ , and $G$ is a finite group acting freely on the product (see \S2 for the full definition). Various rigidity properties of Beauville surfaces have been proved by Catanese~\cite{Cat} and by Bauer, Catanese and Grunewald~\cite{BCG05, BCG06, BCGarX}, including the following (see~\cite[Theorem 4.1]{GT} for a proof by the first and third authors using uniformisation theory):

\begin{proposition}
If $S = (C_1 \times C_2)/G$ and $S' = (C_1' \times C_2') /G'$ are Beauville surfaces with $\pi_1S\cong\pi_1S'$ then $G\cong G'$ and, possibly after transposing factors, each $C_i$ is isomorphic to either $C_i'$ or its complex conjugate curve $\overline {C'_i}$. \hfill$\square$
\end{proposition}

In particular, the conclusions of Proposition~1 apply if $S$ and $S'$ are homeomorphic to each other. Since Beauville surfaces are defined over $\overline{\mathbb{Q}}$, this result (or more precisely its contrapositive) suggests that these surfaces should provide further examples of non-homeomorphic Galois conjugate varieties. Indeed, for this purpose one can use any Beauville surface $S = (C_1 \times C_2)/G$ where $g_1 \ne g_2$ and there is some $\sigma\in\Gal\overline{\mathbb{Q}}/\mathbb{Q}$ such that $C_1^{\sigma}$ is not isomorphic to $C_1$ or $\overline C_1$.

In~\cite{Str} Streit developed a method for determining the Galois orbits and fields of definition of certain curves with large automorphism groups, such as the Macbeath-Hurwitz curves. In~\cite{GT}, the first and third authors used generating triples for the groups $G=PSL_2(p)$, where $p$ is prime, together with Streit's method, to construct arbitrarily large Galois orbits of mutually non-homeomorphic Beauville surfaces. Specifically, for any integer $n>6$ dividing $(p\pm 1)/2$, they constructed an orbit of at least $\varphi(n)/2$ such surfaces, where $\varphi$ is Euler's function. Now the most convenient necessary and sufficient condition for two pairs of triples in a group $G$ to give isomorphic Beauville surfaces depends on a rather delicate relationship between the actions of inner and outer automorphisms of $G$. In this particular case, the existence of a non-trivial outer automorphism of $PSL_2(p)$ (induced by conjugation in $PGL_2(p)$) makes it difficult to determine the precise size of this orbit. Here we use a similar construction, based instead on the groups $G=PGL_2(p)$ which have no outer automorphisms, to give exact values for the (unbounded) sizes of certain Galois orbits of mutually non-homeomorphic Beauville surfaces (see Theorem~2 in \S6). As a particular case, we have following result:

\begin{theorem}\label{thm:particular}
For each prime $p\equiv 19$ mod~$(24)$, there is an orbit of $\Gal\,\overline{\mathbb{Q}}/\mathbb{Q}$ consisting of  $\varphi(m)/4$ Beauville surfaces with Beauville group $PGL_2(p)$, where $m=(p^2-1)/2$, and with mutually non-isomorphic fundamental groups. \hfill$\square$
\end{theorem}

Of course, the surfaces in each such orbit are mutually non-homeomorphic. By Dirichlet's Theorem, there are infinitely many primes $p\equiv 19$ mod~$(24)$. The corresponding orbit-lengths $\varphi(m)/4$ are unbounded above, since there are only finitely many integers $m$ with a given value of $\varphi(m)$.

It is worth noting that the mutually non-isomorphic fundamental groups in Theorems~\ref{thm:particular} and~\ref{thm:mainthm} have isomorphic profinite completions (cf.~Serre's examples in~\cite{Ser64}). We recall that the profinite completion of a group $\Gamma$ is the projective limit $\widehat{\Gamma} = \varprojlim{\Gamma/\Gamma_{i}}$, where $\Gamma_{i}$ ranges over the finite index normal subgroups of $\Gamma$ and the quotient groups are endowed with the obvious epimorphisms $\Gamma/\Gamma_{i}\to \Gamma/\Gamma_{j}$ whenever $\Gamma_{i} <\Gamma_{j}$. In the particular case in which $\Gamma =\pi_{1}X$ is the fundamental group of a non-singular complex projective variety $X$, elementary covering space theory shows that, if $\widetilde{X}$ denotes the holomorphic universal cover of $X$, then $\Gamma$ acts freely and properly discontinuously on $\widetilde{X}$, so that $X$ can be viewed as $X=\widetilde{X}/\Gamma$ and its unramified finite Galois coverings as $X_{i}=\widetilde{X}/\Gamma_{i}$. Therefore $\widehat{\Gamma}$ can be seen as
    \[\widehat{\Gamma}=\widehat{\pi_{1}X}=\varprojlim\Aut(X_{i}/X),\]
where $X_{i}\to X$ ranges over the unramified finite Galois coverings of $X$ by (necessarily) projective varieties $X_{i}$ and the covering groups  $\Aut(X_{i}/X)$ are understood to be related by the natural epimorphisms $\Aut(X_{i}/X)\to\Aut(X_{j}/X)$, whenever $X_{i}\to X$ factors through the coverings $X_{i}\to X_{j}$ and $X_{j}\to X$. Thus, for any Galois element $\sigma$ one clearly has
    \[\widehat{\pi_1X} =\varprojlim\Aut(X_{i}/X)\cong\varprojlim\Aut(X^{\sigma}_{i}/X^{\sigma}) = \widehat{\pi_{1}X^{\sigma}}\]

In \S2 and \S3 we summarise some background information on Beauville surfaces and the groups $G=PGL_2(p)$. In \S4 and \S5 we define and enumerate two different types of generating triples for $G$, which are used in \S6 to construct Beauville surfaces. In \S7 we determine the Galois orbits on these surfaces. This section includes a more general theorem of the above type, along with some further applications and open problems.

\medskip

\noindent{\bf Acknowledgement} The second author is grateful to the Departamento de Matem\'aticas, Universidad Aut\'onoma de Madrid, for financially supporting a visit during which much of this research was carried out.

The first and third authors are partially supported by the grant MTM2009-11848 of the MICINN.

Finally, the first two authors are grateful to the ICMS of Edinburgh, where they first began working on Beauville surfaces.

\section{Beauville surfaces and Beauville structures}\label{sec:beauville}

Beauville surfaces were introduced by Catanese in~\cite{Cat} following an example of Beauville in~\cite{Bea}, and their properties have subsequently been investigated by himself, Bauer and Grunewald~\cite{BCG05, BCG06, BCGarX}. A {\sl Beauville surface\/} (of unmixed type) is a compact complex surface $S$ such that
\begin{description}
\item[\rm(a)] $S$ is isogenous to a higher product, that is, $S\cong (C_1\times C_2)/G$ where $C_1$ and $C_2$ are projective curves of genus at least $2$ and $G$ is a finite group acting freely by holomorphic transformations on $C_1\times C_2$;
\item[\rm(b)] $G$ acts faithfully on each $C_i$ so that $C_i/G$ is isomorphic to the projective line $\mathbb{P}^1(\mathbb{C})$ and the covering $C_i\to C_i/G$ is ramified over three points.
\end{description}

A rational function $C_i\to\mathbb{P}^1(\mathbb{C})$ ramified over at most three points is known as a Bely\u\i\/ function. By Bely\u\i's Theorem~\cite{Bel}, the existence of such a function is equivalent to $C_i$ being defined over $\overline{\mathbb{Q}}$. When, as in condition~(b), a Bely\u\i\/ function is a regular covering, $C_i$ is called a {\sl quasiplatonic curve}.

A group $G$ arises in the above way if and only if it has generating triples $(a_i, b_i, c_i)$ for $i=1, 2$, of orders $(l_i, m_i, n_i)$, such that
\begin{description}
\item[\rm(1)] $a_ib_ic_i=1$ for each $i=1, 2$,
\item[\rm(2)] $l_i^{-1}+m_i^{-1}+n_i^{-1}<1$ for each $i=1, 2$, and
\item[\rm(3)] no non-identity power of $a_1, b_1$ or $c_1$ is conjugate in $G$ to a power of $a_2, b_2$ or $c_2$.
\end{description}

\noindent Such a pair of triples $(a_i, b_i, c_i)$ is called a {\sl Beauville structure\/} on $G$, of bitype $(l_1,m_1,n_1;l_2,m_2,n_2)$. Property~(1) is equivalent to condition~(b), with $a_i$, $b_i$ and $c_i$ representing the local monodromies over the three ramification points. Property~(2) is equivalent to each $C_i$ having genus at least $2$, arising as a smooth quotient $\mathbb{H}/M_i$ of the hyperbolic plane $\mathbb{H}$, where $M_i$ is the kernel of the natural epimorphism $\rho_{i}$ from the triangle group $\Delta_i$ of type $(l_i,m_i,n_i)$ onto $G$. Property~(3) is equivalent to $G$ acting freely on $C_1\times C_2$. It is shown in~\cite{BCG05} that properties~(1) and (3) imply (2).

The fundamental group $\pi_1S$ of a Beauville surface $S$ is the inverse image of the diagonal subgroup under the natural epimorphism $\rho_{1}\times\rho_{2}:\Delta_1\times\Delta_2\to G\times G$, so that
    \[\pi_{1}S\cong\{(\gamma_{1},\gamma_{2})\in\Delta_{1}\times\Delta_{2}:\rho_{1}(\gamma_{1})=\rho_{2}(\gamma_{2})\}\]
It has a normal subgroup $M_1\times M_2\cong \pi_1C_1\times\pi_1C_2\cong\Pi_{g_1}\times\Pi_{g_2}$ with quotient group $G$, where $\Pi_g$ denotes a surface group of genus $g$.

\section{Properties of $PGL_2(p)$}

From now onwards we let
\[G := PGL_2(p)=GL_2(p)/\{\lambda I\mid \lambda\in\F_p^*\},\]
a group of order $p(p^2-1)$, for some prime $p$. This group is complete, i.e.~the centre $Z(G)$ and the outer automorphism group $\Out G$ are both trivial, so $\Aut G\cong G$, acting by conjugation. (See~\cite[\S\S II.6--II.8]{Hup} for properties of $G$.)

Let $p>2$, so that there are three conjugacy classes of maximal cyclic subgroups of $G$. These are
\begin{itemize}
\item elliptic subgroups, of order $p+1$, acting regularly on the projective line $\mathbb{P}^1(p)$ over $\F_p$;
\item parabolic subgroups, of order $p$, with one fixed point and one regular orbit;
\item hyperbolic subgroups, of order $p-1$, with two fixed points and one regular orbit.
\end{itemize}
The elliptic and hyperbolic cyclic subgroups $C$ of order $p\pm 1$ have dihedral normalisers in $G$, of order $2(p\pm 1)$; each element $g\in C$ is conjugate in $G$ to $g^{\pm 1}$, but to no other elements of $C$. The parabolic cyclic subgroups $C$ of order $p$ have normalisers of order $p(p-1)$; these are the stabilisers in $G$ of points in $\mathbb{P}^1(p)$, isomorphic to the affine general linear group $AGL_1(p)$; in this case, all non-identity elements of $C$ are conjugate to each other. In all cases except the involutions, the centraliser in $G$ of a non-identity element is the unique maximal cyclic subgroup containing it; in the case of the involutions, it is the normaliser of that maximal cyclic subgroup, namely a dihedral group
containing it as a subgroup of index 2.

The parabolic elements all lie in the subgroup $G^{+} := PSL_2(p)$ of index $2$ in $G$, whereas elliptic and hyperbolic elements, of order $m$ dividing $p\pm 1$, lie in $G^{+}$ if and only if $(p\pm 1)/m$ is even. It follows that there are two conjugacy classes of involutions in $G$, one of them contained in $G^{+}$ and the other in $G\setminus G^{+}$. Any generating triple for $G$ must contain one element of $G^{+}$, and two of $G\setminus G^{+}$.

For the rest of this paper we let $p \equiv 19$ mod~$(24)$, or equivalently $p\equiv 3$ mod~$(8)$ and $p\equiv 1$ mod~$(3)$.

\section{The first triples}

Let $k$ be any divisor of $p-1$ such that $(p-1)/k$ is odd, so the elements of this order in $G$ all lie in $G\setminus G^{+}$. In particular, since $p \equiv 19$ mod~$(24)$ we can write $k=2k_0$ for some odd number $k_0$.

For each such $k$ let $\T_k$ be the set of all triples $(a_1,b_1,c_1)$ of type $(2,3,k)$ in $G$. Since all elements of odd order lie in $G^{+}$, it follows that each such triple has $a_1\in G\setminus G^{+}$, $b_1\in G^{+}$, and $c_1\in G\setminus G^{+}$. By our choice of $p$, all three elements of such a triple are hyperbolic, and hence so are all their non-identity powers.

\begin{lemma}\label{lem:bc2generateS}
If $(a_1,b_1,c_1)\in\T_k$ with $k>10$ then $b_1$ and $c_1^2$ generate $G^{+}$.
\end{lemma}
\noindent{\sl Proof.}
First let us note that $c_{1}^{2}$ lies in $G^{+}$. It is sufficient to show that no maximal subgroup of $G^{+}$ contains both of these elements. Dickson~\cite[Ch.~XII]{Dic} classified the maximal subgroups of the groups $PSL_2(q)$ for all prime powers $q$. If $q$ is an odd prime $p$ then each maximal subgroup is of one of the following types:
\begin{itemize}
\item the stabiliser of a point in $\mathbb{P}^1(p)$, of order $p(p-1)/2$;
\item a dihedral group of order $p\pm 1$;
\item a subgroup isomorphic to $A_4$, $S_4$ or $A_5$.
\end{itemize}
Both $b_1$ and $c_1^2$ lie in point-stabilisers in $G^{+}$, but not in the same one, for otherwise $b_1$ and $c_1$ would lie in the same point-stabiliser in $G$, isomorphic to $AGL_1(p)$, whereas this group contains no triples of type $(2,3,k)$ for $k>6$. The same argument deals with dihedral subgroups of order $p-1$, except that we replace point-stabilisers in $G$ with dihedral subgroups of order $2(p-1)$. Dihedral subgroups of order $p+1$ are excluded since they have no elements of order $3$, while $A_4$, $S_4$ and $A_5$ have none of order $k/2$ for $k>10$. \hfill$\square$

\begin{corollary}
If $k>10$ then each triple $(a_1,b_1,c_1)\in\T_k$ generates $G$.
\end{corollary}

\noindent{\sl Proof.} This follows immediately from Lemma~\ref{lem:bc2generateS}, since $G^{+}$ is a maximal subgroup of $G$ and $a_1\not\in G^{+}$. \hfill$\square$

\medskip

From now on, we will always assume that $k>10$. There is a natural action of $\Aut G$ on $\T_k$. Since $\Aut G=\Inn G$, this action preserves the conjugacy classes containing the elements of each triple.  By Corollary~1, only the identity automorphism can fix a triple in $\T_k$, so $\Aut G$ acts semiregularly (i.e. freely) on $\T_k$, with $n_k$ orbits where $|\T_k|=n_k|\negthinspace\Aut G|=n_kp(p^2-1)$.

The triples $(a_1,b_1,c_1)\in\T_k$ all have their elements $a_1$ of order $2$ in the same conjugacy class, namely the unique conjugacy class $\A$ of involutions in $G\setminus G^{+}$, and similarly their elements $b_1$ all lie in the unique class $\B$ of elements of order $3$ in $G$. There are $\varphi(k)/2$ conjugacy classes $\C$ of elements of order $k$ in $G$, so for each such class $\C$ let $\T_k(\C)$ denote the set of triples in $\T_k$ with $c_1\in\C$. Thus $\T_k$ is the disjoint union of the sets $\T_k(\C)$, each of which is invariant under $\Aut G$ and is therefore a union of orbits of $\Aut G$.

\begin{lemma}\label{lem:T_k(C)}
For each conjugacy class $\C$ of elements of order $k$ in $G$ we have $|\T_k(\C)|=p(p^2-1)$. \hfill$\square$
\end{lemma}

\noindent{\sl Proof.} We can represent elements of $G$ by pairs $\pm A$ of $2\times 2$ matrices of determinant $\pm 1$ over $\F_p$ (note that $-1$ is not a square in $\F_p$, since $p\equiv 3$ mod~$(4)$). If $(a_1,b_1,c_1)\in\T_k(\C)$ then there are $|\A|=p(p+1)/2$ possible choices for the involution $a_1$, and conjugating by a suitable element of $G$, we can assume that it is represented by the matrix
\[A_1=\Big(\,\begin{matrix}1&0\cr 0&-1\end{matrix}\,\Big).\]
The element $b_1$ of order $3$ is represented by a matrix
\[B_1=\Big(\,\begin{matrix}a&b\cr c&d\end{matrix}\,\Big)\]
with $ad-bc=1$ and (multiplying by $-1$ if necessary) $a+d=1$. Then
\[A_1B_1=\Big(\,\begin{matrix}a&b\cr -c&-d\end{matrix}\,\Big),\]
so $a-d=\pm t$, the trace of a matrix representing elements of the (inverse-closed) class $\C$. Thus
\[a=\frac{1\pm t}{2}\quad{\rm and}\quad d=\frac{1\mp t}{2},\]
so
\begin{equation}
bc=ad-1=\frac{-3-t^2}{4}.
\end{equation}
Now $t^2\ne -3$, for otherwise $bc=0$ and hence $a_1$ and $b_1$ have a common fixed point in $\mathbb{P}^1(p)$, contradicting Corollary~1. It follows that there are $p-1$ solutions $b, c\in\F_p^*$ of equation~(1), and hence (allowing for the two choices for the $\pm$ sign) there are $2(p-1)$ possible elements $b_1$ represented by matrices $B_1$. Multiplying this by the number $p(p+1)/2$ choices for $a_1$, we see that there are $p(p^2-1)$ triples in $\T_k(\C)$. \hfill$\square$

\medskip

Since $\Aut G$ has order $p(p^2-1)$, Lemma~\ref{lem:T_k(C)} shows that each $\T_k(\C)$ is an orbit of this group, so we have:

\begin{corollary}\label{cor:orbits(2,3,k)}
If $k>10$ then $\Aut G$ has $\varphi(k)/2$ orbits on $\T_k$, namely the sets $\T_k(\C)$ where $\C$ ranges over the conjugacy classes of elements of order $k$ in $G$. In particular, the orbit of a triple is characterized by the conjugacy class of its element of order $k$. \hfill$\square$
\end{corollary}

\begin{corollary}\label{cor:orbits(2,3,k)2}
For $k>10$ and a fixed element $c$ of order $k$, we can take representatives of the $\varphi(k)/2$ orbits of $\Aut G$ on $\T_k$ of the form $(a_{i},b_{i},c^{r_{i}})$, for $i=1,\ldots,\varphi(k)/2$, $1\le r_{i} \le k/2$ and $r_{i}$ coprime to $k$.  \hfill$\square$
\end{corollary}

These $\varphi(k)/2$ orbits correspond to the torsion-free normal subgroups $M_i\;(i=1,\ldots, \varphi(k)/2)$ of the triangle group $\Delta_1$ of type $(2,3,k)$ such that $\Delta_1/M_i\cong G$ where, as noted in section~\ref{sec:beauville}, $M_{i}$ is the kernel of the obvious epimorphism $\rho_{i}:\Delta_{1}\to G$ determined by any triple of the corresponding orbit. Let $X_i$ denote the quasiplatonic curve $\mathbb{H}/M_i$ uniformised by $M_i$.

\begin{proposition}\label{prop:(2,3,k)}
The $\varphi(k)/2$ curves $X_i$ have the following properties:
\begin{enumerate}
\item they are mutually non-isomorphic;
\item they all have automorphism group $\Aut X_i\cong G$;
\item they all have the real subfield $K=\mathbb{Q}(\zeta_k)\cap\mathbb{R}$ of the $k$-th cyclotomic field $\mathbb{Q}(\zeta_k)$ as their moduli field and field of definition, where $\zeta_k:=\exp(2\pi i/k)$;
\item they form a single orbit under the Galois group $\Gal K/\mathbb{Q}$.\hfill$\square$
\end{enumerate}
\end{proposition}

We recall that the field of moduli of an algebraic variety $V$ defined over $\overline{\mathbb{Q}}$ is the subfield of $\overline{\mathbb{Q}}$ consisting of all elements fixed by the inertia group $I(V)=\{\sigma\in \mathrm{Gal}(\overline{\mathbb{Q}}/\mathbb{Q}) \mid V^{\sigma}\cong V\}$. The field of moduli is contained in any field of definition, but in general these two fields are not equal. However, Wolfart~\cite{Wol} has shown that quasiplatonic curves are always definable over their fields of moduli.

Let us stress here that in~\cite[Theorem~3]{Str} Streit proves the corresponding results for curves uniformised by normal subgroups of $\Delta_1$ with quotient group isomorphic to $G^{+}=PSL_2(p)$, where $k$ divides $p\pm 1$. His method involves representing curves by their canonical models, and then studying the effect of Galois conjugation on local multipliers, the factors by which automorphisms multiply local coordinates near their fixed points. In order to prove the proposition we will need the following result from~\cite{GT}, which sums up Streit's method in a more general context.

\begin{lemma}\label{lem:galois}Let $G$ be a finite group and $(a,b,c)$ a triple of generators of type $(l,m,n)$ defining a curve $C$. Then for any $\sigma\in\Gal \overline{\mathbb{Q}}/\mathbb{Q}$ the curve $C^{\sigma}$ corresponds to a hyperbolic triple of generators $(a_{\sigma},b_{\sigma},c_{\sigma})$ of $G$ of the form
    \[ a_{\sigma} = ga^{\alpha}g^{-1},\qquad b_{\sigma} = hb^{\beta}h^{-1},\qquad c_{\sigma} = c^{\gamma} \]
where $\sigma(\zeta_{l}^{\alpha})=\zeta_{l}$, $\sigma(\zeta_{m}^{\beta})=\zeta_{m}$ and $\sigma(\zeta_{n}^{\gamma})=\zeta_{n}$ and $g,h\in G$.

In the particular case in which $\sigma$ is complex conjugation, $(a_{\sigma},b_{\sigma},c_{\sigma})=(a^{-1},ab^{-1}a^{-1},c^{-1})$.
\end{lemma}

\noindent{\sl Proof of Proposition~\ref{prop:(2,3,k)}.} (1) We have $X_i\cong X_{j}$ if and only if $M_i^\gamma=M_{j}$ for some $\gamma\in PSL_2(\mathbb{R})$. If this is the case then $N(M_i)^\gamma=N(M_{j})$, where $N(\;)$ denotes the normaliser in $PSL_2(\mathbb{R})$. Now $M_i$ is normal in $\Delta_1$, so $N(M_i)$ is a Fuchsian group containing $\Delta_1$. By Singerman's classification~\cite{Sin}, the triangle group of type $(2,3,k)$ is a maximal Fuchsian group for $k>6$, so $N(M_i)=\Delta_1$, and similarly $N(M_{j})=\Delta_1$. Thus $\Delta_1^\gamma=\Delta_1$, so $\gamma\in N(\Delta_1)=\Delta_1$ and hence $M_i=M_{j}$, giving $i=j$.

(2) We have $\Aut X_i\cong N(M_i)/M_i$. The argument used to prove (1) shows that $N(M_i)=\Delta_1$, so $\Aut X_i\cong \Delta_1/M_i\cong G$.

(3) Let the triple $(a_{1},b_{1},c)\in\T_{k}(\C)$ correspond to $X_{1}$. In view of the definition of a field of moduli, the first part of Lemma~\ref{lem:galois} clearly implies that the moduli field of $X_{1}$ is contained in $\mathbb{Q}(\zeta_k)$, and the second part of it states that the complex conjugate curve $\overline{X_{1}}$ is defined by $(a^{-1},ab^{-1}a^{-1},c^{-1})$, which lies in $\T_{k}(\C)$ too, and therefore $\overline{X_{1}}\cong X_{1}$. As a consequence the moduli field of $X_{1}$ is contained in $K=\mathbb{Q}(\zeta_k)\cap\mathbb{R}$.\\ On the other hand, for every triple $(a',b',c')\in\T_{k}$, defining a curve $X'$, we can suppose, by Corollary~\ref{cor:orbits(2,3,k)2}, that $c'=c^{r}$ for some $r$ coprime to $k$. Now, by Lemma~\ref{lem:galois}, for any element $\sigma\in\Gal \overline{\mathbb{Q}}/\mathbb{Q}$ such that $\sigma(\zeta_{k}^{r})=\zeta_{k}$ one has $(a_{\sigma},b_{\sigma},c_{\sigma})=(ga^{\alpha}g^{-1}, hb^{\beta}h^{-1},c^{r})$, and by Corollary~\ref{cor:orbits(2,3,k)} it follows that $(a_{\sigma},b_{\sigma},c_{\sigma})$ and $(a',b',c')$ are $\mathrm{Aut}\,G$-equivalent. Hence $X_{i}^{\sigma}=X'$ and as a consequence the $\varphi(k)/2$ curves $X_{1},\ldots,X_{\varphi(k)/2}$ are Galois conjugate. \\ Now, let us note that the field of moduli of a quasiplatonic curve is always a field of definition of such a curve (see~\cite{Wol}), and therefore its degree is always greater than or equal to the cardinality of the Galois orbit of $X_{1}$. Since the field $K$ has exactly degree $\varphi(k)/2$, it follows that $K$ is the field of moduli (hence field of definition) of $X_{1}$, and therefore of each $X_{i}$.

(4) This follows from the proof of (3).
\hfill$\square$

{\subsection{An alternative proof of Lemma~\ref{lem:T_k(C)}}

Here we outline an alternative method of proof of Lemma~2 using character theory, which may be useful in groups where calculations with explicit elements, as above, are not so straightforward (see e.g.~\cite{FJ}). We use the following well-known result (see~\cite[\S7.2]{Ser92} for this and other similar results):

\begin{proposition}
If $\A$, $\B$ and $\C$ are conjugacy classes in a finite group $G$, then the number of solutions $(a,b,c)\in\A\times\B\times\C$ of the equation $abc=1$ is given by the formula
\[\frac{|\A||\B||\C|}{|G|}\sum_{\chi}\frac{\chi(a)\chi(b)\chi(c)}{\chi(1)},\]
where $\chi$ ranges over the irreducible complex characters of $G$. \hfill$\square$
\end{proposition}

The character table for $G=PGL_2(p)$ can be obtained from the generic character table for $GL_2(q)$ for all prime powers $q$ (see~\cite[\S15.9]{DM}, for instance) by putting $q=p$ and restricting attention to those irreducible characters of $GL_2(q)$ which are constant on the scalar matrices, so that they correspond to representations of $G$.

In the case of Lemma~2 we have $|\A|=p(p+1)/2$, $|\B|=|\C|=p(p+1)$ and $|G|=p(p^2-1)$, so
\[|\T_k(\C)|=
\frac{p^2(p+1)^2}{2(p-1)}\sum_{\chi}\frac{\chi(a_1)\chi(b_1)\chi(c_1)}{\chi(1)}.\]
The character table for $G$ shows that as $p\to\infty$ the sum on the right-hand side is dominated by the two characters $\chi$ of degree $1$ (those of $G/G^{+}\cong C_2$), which each contribute $1$ to the summation. (More precise estimates of the character sum are aided by the fact that nearly half of the characters $\chi$, specifically those of degree $p-1$, take the value $0$ on hyperbolic elements, so they contribute nothing to the sum.)
Thus
\[\frac{|\T_k(\C)|}{|G|}=
\frac{p(p+1)}{2(p-1)^2}\sum_{\chi}\frac{\chi(a_1)\chi(b_1)\chi(c_1)}{\chi(1)}\]
approaches 1 as $p\to\infty$. But this number is an integer, the number of (regular) orbits of $\Aut G=G$ on $T_k(\C)$, so it must be equal to $1$, giving $|\T_k(\C)|=|G|=p(p^2-1)$. This argument provides a proof of Lemma~\ref{lem:T_k(C)} valid for sufficiently large primes $p\equiv 19$ mod~$(24)$, but the careful proof it outlines, using exact character values, is valid for all such $p$.

\section{The second triples}

Now let $l$ be any divisor of $p+1$ such that $(p+1)/l$ is odd. In this case, our choice of $p$ implies that there is an odd number $l_0$ such that $l = 4l_0$. If $(a_2,b_2,c_2)$ is any triple of type $(2,4,l)$ in $G$, then $a_2\in G^{+}$, $b_2\in G\setminus G^{+}$, and $c_2\in G\setminus G^{+}$. In this case the non-identity powers of $a_2$, $b_2$ and $c_2$ are all elliptic. Arguments similar to those used in the preceding section show that provided $l>10$, each such triple generates $G$ and there are $\varphi(l)/2$ orbits of $\Aut G$ on such triples, one for each of the $\varphi(l)/2$ conjugacy classes of elements $c_2$ of order $l$ in $G$. (The involution $a_2$, represented by the matrix
\[A_2=\Big(\,\begin{matrix}0&1\cr -1&0\end{matrix}\,\Big),\]
has $p(p-1)/2$ conjugates, and the solutions of the analogue of equation~(1) form two quadrics, each with $p+1$ points.) Since $l>8$ the triangle group $\Delta_2$ of type $(2, 4, l)$ is a maximal Fuchsian group~\cite{Sin}. Replacing the generator $b_1$ of order $3$ with $b_2$ of order $4$ is not significant,  so as in the case of the first triples we find that the quasiplatonic curves $Y_j$ corresponding to these orbits of triples satisfy:

\begin{proposition}
The $\varphi(l)/2$ curves $Y_j$ have the following properties:
\begin{enumerate}
\item they are mutually non-isomorphic;
\item they all have automorphism group $\Aut Y_j\cong G$;
\item they all have the real subfield $L=\mathbb{Q}(\zeta_l)\cap\mathbb{R}$ of the $l$-th cyclotomic field $\mathbb{Q}(\zeta_l)$ as their moduli field and field of definition;
\item they form a single orbit under the Galois group $\Gal L/\mathbb{Q}$ of $L$. \hfill$\square$
\end{enumerate}
\end{proposition}

One can also apply the alternative argument given in \S4.1, with minor modifications, to the triples of type $(2,4,l)$ considered here: in this case the characters of degree $p+1$ vanish on the elliptic elements. This type of argument explains why we needed to choose both of the generating triples in $G$ to include involutions: otherwise, we would have $|\A|, |\B|, |\C|\sim p^2$ as $p\to\infty$ and hence $\Aut G$ would have two orbits, rather than one, on generating triples in $\A\times\B\times\C$.

\section{The Beauville surfaces}

If $(a_1,b_1,c_1)$ and $(a_2,b_2,c_2)$ are triples in $G$ of types $(2,3,k)$ and $(2,4,l)$, with $k, l>10$, then since the non-identity powers of $a_1, b_1$ and $c_1$ are hyperbolic, while those of $a_2, b_2$ and $c_2$ are elliptic, these two triples form a Beauville structure of bitype $(2,3,k;2,4,l)$, corresponding to a Beauville surface
\[S_{ij}=(X_i\times Y_j)/G.\]
Since $k = 2k_0$ and $l = 4l_0$ for coprime odd $k_0$ and $l_0$, the number of such surfaces $S_{ij}$ is
\[\frac{\varphi(k)}{2}.\frac{\varphi(l)}{2}
=\frac{\varphi(k_0)}{2}.\frac{2\varphi(l_0)}{2}
= \frac{\varphi(k_0l_0)}{2}
= \frac{\varphi(m)}{4},\]
where
\[m = {\rm lcm}(k,l) = 4k_0l_0.\]
By Proposition~2 the $\varphi(k)/2$ curves $X_i$ are real and mutually non-isomorphic, as are the $\varphi(l)/2$ curves $Y_j$ by Proposition~3. No pair $X_i$ and $Y_j$ can be isomorphic, since they are uniformised by surface groups with non-isomorphic normalisers $\Delta_1$ and $\Delta_2$. It therefore follows from Proposition~1 that the $\varphi(m)/4$ surfaces $S_{ij}$ have mutually non-isomorphic fundamental groups. In particular, they are mutually non-homeomorphic.

Moreover, up to isomorphism there cannot be any more Beauville surfaces with group $G$ and bitype $(2,3,k;2,4,l)$. This is because if there was another Beauville surface $S'$, its defining triples $(a_{1}',b_{1}',c_{1}')$ and $(a_{2}',b_{2}',c_{2}')$ would be conjugate to the two triples defining one of our surfaces $S_{ij}$ by means of elements $g_{1},g_{2}\in G$. Now, if for $r=1,2$ we choose a preimage $\gamma_{r}\in\Delta_{r}$ of $g_{r}$ under the epimorphism $\rho_{r}:\Delta_{r}\to G$ determined by the triple $(a_{r}',b_{r}',c_{r}')$ then, clearly the groups $\pi_{1}S'$ and $\pi_{1}S_{ij}$ uniformising the surfaces $S'$ and $S_{ij}$ (see section~\ref{sec:beauville}) are conjugate under the element $(\gamma_{1},\gamma_{2})\in\mathrm{Aut}(\mathbb{H}\times\mathbb{H})$. As a consequence, we can characterize the surface $S_{ij}$ as the only Beauville surface with group $G$, bitype $(2,3,k;2,4,l)$ and curves $X_i$ and $Y_j$.

\begin{example}For $p=19$ we can take $k=18$ and $l=20$. By the results of the previous sections there are $\varphi(18)/2=3$ orbits of $\mathrm{Aut}\,G$ on triples of generators of $G=PGL_{2}(19)$ of type $(2,3,18)$, and $\varphi(20)/2=4$ orbits on triples of type $(2,4,20)$. By computer means we can find representatives
\begin{eqnarray*}
  (a_{1},b_{1},c_{1}) &=& \left(
  \left( \begin{array}{cc}  6 & 12 \\ 5 & 13 \\ \end{array} \right),
  \left( \begin{array}{cc}  3 & 12 \\ 12 & 13 \\ \end{array} \right),
  \left( \begin{array}{cc}  2 & 0 \\ 0 & 1 \\ \end{array} \right)\right) \\
  (a_{1}',b_{1}',c_{1}^{5}) &=& \left(
  \left( \begin{array}{cc}  2 & 6 \\ 9 & 17 \\ \end{array} \right),
  \left( \begin{array}{cc}  6 & 6 \\ 8 & 17 \\ \end{array} \right),
  \left( \begin{array}{cc}  13 & 0 \\ 0 & 1 \\ \end{array} \right)\right) \\
  (a_{1}'',b_{1}'',c_{1}^{7}) &=& \left(
  \left( \begin{array}{cc}  11 & 10 \\ 7 & 8 \\ \end{array} \right),
  \left( \begin{array}{cc}  13 & 10 \\ 10 & 8 \\ \end{array} \right),
  \left( \begin{array}{cc}  14 & 0 \\ 0 & 1 \\ \end{array} \right)\right) \\
\end{eqnarray*}
of the first three orbits, defining curves $X_{1},X_{2},X_{3}$, and representatives
\begin{eqnarray*}
  (a_{2},b_{2},c_{2}) &=& \left(
  \left( \begin{array}{cc}  0 & 9 \\ 2 & 0 \\ \end{array} \right),
  \left( \begin{array}{cc}  10 & 9 \\ 2 & 15 \\ \end{array} \right),
  \left( \begin{array}{cc}  1 & 2 \\ 1 & 1 \\ \end{array} \right)\right) \\
  (a_{2}',b_{2}',c_{2}^{3}) &=& \left(
  \left( \begin{array}{cc}  7 & 11 \\ 11 & 12 \\ \end{array} \right),
  \left( \begin{array}{cc}  13 & 7 \\ 17 & 12 \\ \end{array} \right),
  \left( \begin{array}{cc}  7 & 10 \\ 5 & 7 \\ \end{array} \right)\right) \\
  (a_{2}'',b_{2}'',c_{2}^{7}) &=& \left(
  \left( \begin{array}{cc}  13 & 1 \\ 1 & 6 \\ \end{array} \right),
  \left( \begin{array}{cc}  12 & 6 \\ 4 & 13 \\ \end{array} \right),
  \left( \begin{array}{cc}  11 & 15 \\ 17 & 11 \\ \end{array} \right)\right) \\
  (a_{2}''',b_{2}''',c_{2}^{9}) &=& \left(
  \left( \begin{array}{cc}  11 & 7 \\ 7 & 8 \\ \end{array} \right),
  \left( \begin{array}{cc}  11 & 13 \\ 9 & 14 \\ \end{array} \right),
  \left( \begin{array}{cc}  6 & 13 \\ 16 & 6 \\ \end{array} \right)\right) \\
\end{eqnarray*}
of the last four orbits, defining curves $Y_{1},Y_{2},Y_{3},Y_{4}$. Any other triple $(r,s,t)$ of type $(2,3,18)$ or $(2,4,20)$ can be mapped by an automorphism of $PGL_{2}(19)$ into one of the first three or last four orbits, depending on the conjugacy class of $t$.
Consequently, we can construct $12$ pairwise non-isomorphic Beauville surfaces of the form $S_{ij}=(X_{i}\times Y_{j})/PGL_{2}(19)$, where $1\le i\le 3$ and $1\le j\le 4$.
\end{example}

\section{The Galois orbits}

Since $K\cap L=\mathbb{Q}$, the compositum $M$ of $K$ and $L$, i.e.~the subfield $KL$ of $\overline{\mathbb{Q}}$ which they generate, has degree
\[|M:\mathbb{Q}| = |K:\mathbb{Q}||L:\mathbb{Q}|
= \frac{\varphi(k)}{2}.\frac{\varphi(l)}{2}=\frac{\varphi(m)}{4}\]
over $\mathbb{Q}$. Since $K$ and $L$ are abelian extensions of $\mathbb{Q}$, so is $M$ (it is, in fact, a proper subfield
of $\mathbb{Q}(\zeta_m)\cap\mathbb{R}$). The Galois group $\Gal M/\mathbb{Q}$ of $M$ over $\mathbb{Q}$ therefore has the form
\[\Gal M/\mathbb{Q}=\Gal M/L \times\Gal M/K
\cong \Gal K/\mathbb{Q} \times \Gal L/\mathbb{Q}.\]
Since the direct factors act regularly on the sets of curves $X_i$ and $Y_j$, it follows that $\Gal M/\mathbb{Q}$ acts regularly on the set of surfaces $S_{ij}$. These surfaces therefore form an orbit $\Omega=\Omega(p,k,l)$ of length $\varphi(m)/4$ under the absolute Galois group. We have thus proved:

\begin{theorem}\label{thm:mainthm}
For each prime $p\equiv 19$ mod~$(24)$, and for each pair of divisors $k, l>10$ of $p-1$ and $p+1$ such that $(p-1)/k$ and $(p+1)/l$ are odd, there is an orbit of $\Gal\overline{\mathbb{Q}}/\mathbb{Q}$ consisting of  $\varphi(m)/4$ Beauville surfaces with Beauville group $PGL_2(p)$, where $m={\rm lcm}(k,l)$, and with mutually non-isomorphic fundamental groups. \hfill$\square$
\end{theorem}

For any prime $p\equiv 19$ mod~$(24)$, one can satisfy the hypotheses of Theorem~2 by taking $k=p-1$ and $l=p+1$, thus proving Theorem~1 (see \S 1). The resulting Beauville structures have bitype $(2,3,p-1;2,4,p+1)$, so that different Galois orbits $\Omega(p,p-1,p+1)$ correspond to Beauville structures of different bitypes. For a given $p$, the fundamental groups of these surfaces are all subgroups of index $|G|=p(p^2-1)$ in $\Delta_1\times\Delta_2$; although they are mutually non-isomorphic, each is an extension of $\Pi_g\times\Pi_h$ by $G$, where the curves $X_i$ and $Y_j$ have genera
\[g=\frac{1}{12}(p-1)(p^2-5p-12)\quad{\rm and}\quad h=\frac{1}{8}(p+1)(p^2-5p+8)\]
by the Riemann-Hurwitz formula.

The most general set of triples $(p,k,l)$ satisfying the conditions of Theorem~2 arises as follows. Given any pair $k, l>10$ with $k=2k_0$ and $l=4l_0$ for coprime odd $k_0$ and $l_0$, the latter coprime to $3$, the congruences
\[p\equiv 19 \; {\rm mod}\,(24),\quad p\equiv k+1\;  {\rm mod}\;(2k),\quad p\equiv l-1\; {\rm mod}\,(2l)\]
are equivalent to
\[p\equiv 3\;{\rm mod}\,(8), \quad p\equiv 1\;{\rm mod}\,(k_0'), \quad p\equiv -1\;{\rm mod}\,(l_0)\]
where $k_0'={\rm lcm}(3,k_0)$, and hence (since $8, k_0'$ and $l_0$ are mutually coprime) to a single congruence mod~$(8k_0'l_0)$, satisfied by infinitely many primes $p$. For example, one could fix $k$ and $l$, so that the bitype $(2,3,k;2,4,l)$ and hence the group $\Delta_1\times\Delta_2$ are fixed, by taking $k=18$ and $l=20$ for primes $p\equiv 19$ mod~$(360)$ for instance, but then the size $\varphi(k)\varphi(l)/4$ of the orbits $\Omega(p,k,l)$ is also fixed. This raises the question of whether there exist arbitrarily large Galois orbits of mutually non-homeomorphic Beauville surfaces, all corresponding to Beauville structures of the same bitype.

The kernel of the action of $\Gal\overline{\mathbb{Q}}/\mathbb{Q}$ on a single orbit $\Omega=\Omega(p,k,l)$ is the subgroup $\Gal\overline{\mathbb{Q}}/M$ where $M=KL$. Note that $KL$ is the moduli field of the surfaces $S_{ij}\in\Omega$. This is because $I(S_{ij})=I(X_{i}\times Y_{j})=I(X_{i})\cap I(Y_{j})$ and therefore the field of moduli of $S_{ij}$ contains the compositum $KL$ of the fields of moduli of $X_{i}$ and $Y_{j}$ while, on the other hand, the subfield fixed by $I(X_{i}\times Y_{j})$ is included in $KL$, since it is a field of definition of $X_{i}\times Y_{j}$. The kernel of the action of $\Gal\overline{\mathbb{Q}}/\mathbb{Q}$ on the union of all the orbits $\Omega(p,k,l)$ is therefore the intersection of these subgroups $\Gal\overline{\mathbb{Q}}/M$. This is $\Gal\overline{\mathbb{Q}}/\M$ where $\M$ is the compositum of all the corresponding moduli fields $M$, a proper subfield of the maximal cyclotomic field
\[\mathbb{Q}^{\rm ab}=\bigcup_{n\in\N}\mathbb{Q}(\zeta_n).\]
This raises the problem of determining the kernel of the action of the absolute Galois group on {\em all\/} Beauville surfaces and, in particular, the question originally posed by Catanese in~\cite{BCGarX} of whether  this action is faithful.


\bigskip


\noindent {\sc G. Gonz\'{a}lez-Diez:} \\ Departamento de Matem\'aticas, Universidad Aut\'onoma de Madrid, 28049, Madrid, Spain. \\ email: {\tt gabino.gonzalez@uam.es}

\

\noindent {\sc G. A. Jones:} \\ School of Mathematics, University of Southampton, Southampton SO17 1BJ, U.K. \\ email: {\tt G.A.Jones@maths.soton.ac.uk}

\

\noindent {\sc D. Torres-Teigell:} \\Departamento de Matem\'aticas, Universidad Aut\'onoma de Madrid, 28049, Madrid, Spain. \\ email: {\tt david.torres@uam.es}

\end{document}